\documentstyle[12pt,amssymb]{article}

\begin{document}
\baselineskip15pt

\newtheorem{definition}{Definition $\!\!$}
\newtheorem{prop}[definition]{Proposition $\!\!$}
\newtheorem{lem}[definition]{Lemma $\!\!$}
\newtheorem{corollary}[definition]{Corollary $\!\!$}
\newtheorem{theorem}[definition]{Theorem $\!\!$}
\newtheorem{example}[definition]{Example $\!\!$}
\newtheorem{remark}[definition]{Remark $\!\!$}

\newcommand{\nc}[2]{\newcommand{#1}{#2}}
\newcommand{\rnc}[2]{\renewcommand{#1}{#2}}
\nc{\bpr}{\begin{prop}}
\nc{\bth}{\begin{theorem}}
\nc{\ble}{\begin{lem}}
\nc{\bco}{\begin{corollary}}
\nc{\bre}{\begin{remark}}
\nc{\bex}{\begin{example}}
\nc{\bde}{\begin{definition}}
\nc{\ede}{\end{definition}}
\nc{\epr}{\end{prop}}
\nc{\ethe}{\end{theorem}}
\nc{\ele}{\end{lem}}
\nc{\eco}{\end{corollary}}
\nc{\ere}{\end{remark}}
\nc{\eex}{\end{example}}
\nc{\lall}{{\mbox{\large $\forall$}}}
\nc{\lex}{{\mbox{\large $\exists$}}}
\nc{\epf}{\hfill\mbox{$\Box$}}
\nc{\ot}{\otimes}
\nc{\bsb}{\begin{Sb}}
\nc{\esb}{\end{Sb}}
\nc{\ct}{\mbox{${\cal T}$}}
\nc{\ctb}{\mbox{${\cal T}\sb B$}}
\nc{\bcd}{\[\begin{CD}}
\nc{\ecd}{\end{CD}\]}
\nc{\bea}{\begin{eqnarray*}}
\nc{\eea}{\end{eqnarray*}}
\nc{\be}{\begin{enumerate}}
\nc{\ee}{\end{enumerate}}
\nc{\beq}{\begin{equation}}
\nc{\eeq}{\end{equation}}
\nc{\bi}{\begin{itemize}}
\nc{\ei}{\end{itemize}}
\nc{\kr}{\mbox{Ker}}
\nc{\te}{\!\ot\!}
\nc{\bmlp}{\mbox{\boldmath$\left(\right.$}}
\nc{\bmrp}{\mbox{\boldmath$\left.\right)$}}
\rnc{\phi}{\mbox{$\varphi$}}
\nc{\LAblp}{\mbox{\LARGE\boldmath$($}}
\nc{\LAbrp}{\mbox{\LARGE\boldmath$)$}}
\nc{\Lblp}{\mbox{\Large\boldmath$($}}
\nc{\Lbrp}{\mbox{\Large\boldmath$)$}}
\nc{\lblp}{\mbox{\large\boldmath$($}}
\nc{\lbrp}{\mbox{\large\boldmath$)$}}
\nc{\blp}{\mbox{\boldmath$($}}
\nc{\brp}{\mbox{\boldmath$)$}}
\nc{\LAlp}{\mbox{\LARGE $($}}
\nc{\LArp}{\mbox{\LARGE $)$}}
\nc{\Llp}{\mbox{\Large $($}}
\nc{\Lrp}{\mbox{\Large $)$}}
\nc{\llp}{\mbox{\large $($}}
\nc{\lrp}{\mbox{\large $)$}}
\nc{\lbc}{\mbox{\Large\boldmath$,$}}
\nc{\lc}{\mbox{\Large$,$}}
\nc{\bc}{\mbox{\boldmath$,$}}
\rnc{\epsilon}{\varepsilon}
\rnc{\ker}{\mbox{\em Ker}}
\nc{\ra}{\rightarrow}
\nc{\ci}{\circ}
\nc{\cc}{\!\ci\!}
\nc{\T}{\mbox{\sf T}}
\nc{\can}{\mbox{\em\sf T}\!\sb R}
\nc{\cnl}{$\mbox{\sf T}\!\sb R$}
\nc{\lra}{\longrightarrow}
\nc{\M}{\mbox{Map}}
\rnc{\to}{\mapsto}
\nc{\imp}{\Rightarrow}
\rnc{\iff}{\Leftrightarrow}
\nc{\bmq}{\cite{bmq}}
\nc{\ob}{\mbox{$\Omega\sp{1}\! (\! B)$}}
\nc{\op}{\mbox{$\Omega\sp{1}\! (\! P)$}}
\nc{\oa}{\mbox{$\Omega\sp{1}\! (\! A)$}}
\nc{\inc}{\mbox{$\,\subseteq\;$}}
\nc{\de}{\mbox{$\Delta$}}
\nc{\spp}{\mbox{${\cal S}{\cal P}(P)$}}
\nc{\dr}{\mbox{$\Delta_{R}$}}
\nc{\dsr}{\mbox{$\Delta_{\cal R}$}}
\nc{\m}{\mbox{m}}
\nc{\hsp}{\hspace*}
\nc{\nin}{\mbox{$n\in\{ 0\}\!\cup\!{\Bbb N}$}}
\nc{\ha}{\mbox{$\alpha$}}
\nc{\hb}{\mbox{$\beta$}}
\nc{\hg}{\mbox{$\gamma$}}
\nc{\hd}{\mbox{$\delta$}}
\nc{\he}{\mbox{$\varepsilon$}}
\nc{\hz}{\mbox{$\zeta$}}
\nc{\hs}{\mbox{$\sigma$}}
\nc{\hk}{\mbox{$\kappa$}}
\nc{\hm}{\mbox{$\mu$}}
\nc{\hn}{\mbox{$\nu$}}
\nc{\hl}{\mbox{$\lambda$}}
\nc{\hG}{\mbox{$\Gamma$}}
\nc{\hD}{\mbox{$\Delta$}}
\nc{\th}{\mbox{$\theta$}}
\nc{\Th}{\mbox{$\Theta$}}
\nc{\ho}{\mbox{$\omega$}}
\nc{\hO}{\mbox{$\Omega$}}
\nc{\hp}{\mbox{$\pi$}}
\nc{\hP}{\mbox{$\Pi$}}

\def\esl{{\mbox{$E\sb{{\frak s}{\frak l} (2,{\Bbb R})}$}}}
\def\esu{{\mbox{$E\sb{{\frak s}{\frak u}(2)}$}}}
\def\zf{{\mbox{${\Bbb Z}\sb 4$}}}
\def\zt{{\mbox{$2{\Bbb Z}\sb 2$}}}

\title{\vspace*{-30mm}\large\bf
The Einstein Action for  Algebras of Matrix Valued Functions --- Toy Models
}
\author{\normalsize Piotr M.~Hajac \thanks{\small
On leave from: Department of Mathematical Methods in Physics,
Warsaw University, ul.~Ho\.{z}a 74, Warsaw, 00--682 Poland.
http://info.fuw.edu.pl/KMMF/ludzie\underline{~~}ang.html 
(E-mail: pmh@fuw.edu.pl)
}\\
\normalsize Mathematics Section, International Centre for Theoretical Physics,\\
\normalsize  Strada Costiera 11, 34014 Trieste, Italy.
}
\date{\normalsize 3 April 1996}
\maketitle

\begin{abstract}
Two toy models are considered within the framework of noncommutative
differential
geometry. In the first one, the Einstein action of the Levi--Civita connection
is computed for the algebra of matrix valued functions on a torus.
It is shown that, assuming some constraints on the metric,
this action splits into a classical-like, a quantum-like and a mixed term.
In the second model, an analogue of the Palatini method of variation is applied
to obtain critical points of the Einstein action functional for
$M\sb 4({\Bbb R})$. It is pointed out that a solution to the Palatini
variational problem is not necessarily a Levi--Civita connection.
In this model, no additional assumptions regarding metrics are made.
\end{abstract}


\section*{I. Introduction}

The goal of this note is to analyse the behaviour of a noncommutative analogue
of the Einstein--Hilbert action functional on two toy models.
General definitions and constructions employed to study those models
are provided in the next section.

In Section~III, we present some results regarding the computation
of the Einstein action of the Levi--Civita connection
for \mbox{$C\sp{\infty}(T\sp m)\te
M\sb n({\Bbb R})\,$}, the algebra of matrix valued functions on an $m$-tori.
The approach proposed there is analogous to the `derivation based' approach to
the calculation of the Yang--Mills (Maxwell) action for
an algebra of matrix valued functions that was
carried out in \cite{dkm} (see Section~V in~\cite{dkm}, cf.~\cite{dvkm}
and the sections 4 and 5 of~\cite{md}).
We choose our manifold $M$ to be an $m$-tori because it is a compact Abelian
group,
and we want integrals over $M$ to be finite, and
$\mbox{Der}( C\sp\infty(T\sp m))$ to be commutative as a Lie algebra
and free as a $C\sp\infty(T\sp m)$-module.
 We also assume that the metric, 
understood as a 
pairing of derivations, has its values in the centre of an algebra.
Consequently, the results presented there can be interpreted as
concerning the `commutative part' of noncommutative geometry. (Indeed,
the Levi--Civita connections for $M\sb n({\Bbb R})$ 
can be interpreted as the torsion part of the flat
connection on $SL(n,{\Bbb R})$ given by the left translations; see
Remark~\ref{flat}).

In Section~IV, which is the main part of this work, we study a toy model that is
based on the algebra $M\sb 4({{\Bbb R}})$ of 4 by 4 real matrices, and the
2-dimensional Lie
subalgebra \mbox{${\frak s}{\frak o}(2)\oplus{\frak s}{\frak o}(2)$} of
$\mbox{Der}(M\sb 4({\Bbb R}))$.
It is a simple but quite computable model that is presented with
the aim of providing hints on how to approach more complicated situations.
For this model, we derive an analogue of the Einstein vacuum field
equation, and apply the Palatini method of variation to obtain critical points.
We find a solution to the Palatini variational problem that is,
in general, neither metric nor torsion free. Yet, this solution (a connection
$\nabla\sp g$ determined by a metric~$g$) turns out to be Ricci flat for all
metrics. Consequently, the Einstein field equation, which depends essentially
on the Ricci curvature (see (\ref{field})), is automatically satisfied, and
$(g,\nabla\sp g)$ is a critical point of the Einstein action functional for
any metric $g$. The value of this functional at any such critical point
is the same (zero). 
Thus we obtain a result that partially
 reflects the classical geometric phenomenon that,
for any Riemannian 2-manifold, the Einstein--Hilbert action computed for a
metric $g$
and its Levi--Civita connection does not depend on $g$ (see 9.1.10 in
\cite{db}).
The key difference between Section~IV and the preceding one is that in
Section~IV we no
longer assume that the metric is centre valued.

The proofs or calculations
are rather straightforward and are often omitted here for the sake of brevity.
Except for Proposition~\ref{spera}, the letter
$k$ will denote a field, $A$ a unital associative $k$-algebra,
$Z(A)$ its centre, and $\hO(\! A)$ a differential graded algebra
\mbox{$A\oplus\bigoplus\sb{n\geq 1}\sp\infty
\mbox{Hom}\sb{Z(A)}(\mbox{\large$\wedge$}\sp n{\cal L},\, A)$}
with the differential defined as in Proposition~\ref{spera}.
The Einstein convention of summing over repeating indices is assumed.

\section*{II. Preliminaries}
{\parindent=0pt

Let ${\cal L}$ be a $Z(A)$-submodule and Lie subalgebra of $\mbox{Der}(A)$ such
that
${\cal L}\ot\sb{Z(A)}A$ is a finitely generated projective right $A$-module.
(Compare with the notion of a Lie--Cartan pair introduced in~\cite{kastora}.)

\bde\label{lcdef}
A linear map
$\nabla : {\cal L}\ot\sb{Z(A)}\Omega\sp *\! (\! A)\lra
{\cal L}\ot\sb{Z(A)}\Omega\sp{*+1}\! (\! A)$ is called a connection 
on~${\cal L}$~iff
\[
\lall\; X\in{\cal L},\, \ha\in\Omega(\! A)\, :\;
\nabla(X\ot\sb{Z(A)}\ha )=(\nabla X )\ha +X\ot\sb{Z(A)}d\ha\, .
\]
\ede

\bde
The endomorphism
\mbox{$\nabla\sp 2\!\in\!
\mbox{\em End}\sb{\Omega(A)}\!\llp\!{\cal L}\te\sb{Z(A)}\hO (\! A)\!\lrp$}
is called the curvature of a connection~$\nabla$.
\ede

\bre\label{usual}\em
The notions of connection and curvature defined above are equivalent to the
usual notions of noncommutative connection and curvature on a projective module
(e.g., see Section~III.B of \cite{dkm}).
In this case, the projective right $A$-module is ${\cal L}\te\sb{Z(A)}A$.
\hfill{$\Diamond$}\ere

\bpr[cf.~p.369 in \cite{ms}]\label{spera}
Let $A$ be an associative unital algebra over a commutative ring $k$.
Let ${\cal L}\sb k$ be a $k$-Lie subalgebra of the space of all $k$-derivations
of
$A$, and let $\cal E$ be any right $A$-module admitting a connection
(see Remark~\ref{usual}). If $\hO(\! A)$ is a differential graded subalgebra of
\mbox{$A\oplus\bigoplus\sb{n\geq 1}\sp\infty                         
\mbox{\em Hom}\sb k(\mbox{\large$\wedge$}\sp n{\cal L}\sb k,\, A)$}
with the differential
given by (see the first section in~\cite{dvkm}):
\bea
(d\ha)(X\sb 0,X\sb 1,\cdots ,X\sb n)\!\!\!
&=&\!\!\!\!\sum\sb{0\leq i\leq n}(-)\sp iX\sb i\ha
(X\sb 0,\cdots,X\sb{i-1},X\sb{i+1},\cdots ,X\sb n)\\
&+&\!\!\!\!\!\!\sum\sb{0\leq r<s\leq n}\!\!\! (-)\sp{r+s}\ha
([X\sb r,X\sb s],X\sb 0,\cdots,X\sb{r-1},X\sb{r+1},
\cdots ,X\sb{s-1},X\sb{s+1},\cdots ,X\sb n)\, ,
\eea
then
\[
\lall\,\xi\in{\cal E},\; X,Y\in{\cal L}\sb k:\;\;
(\nabla\sp 2\xi)(X,Y)=
\left([\nabla\sb X,\nabla\sb Y]-\nabla\sb{[X,Y]}\right)(\xi)\, ,
\]
where, as in the classical differential geometry,
 $\nabla\!\sb Z\,\xi$ denotes $(\nabla\xi)(Z)$.
\epr

\bde\label{metric}
A map $g: ({\cal L}\te\sb{Z(A)}A)\times({\cal L}\te\sb{Z(A)}A)\lra A$
is called a pseudo-Riemannian metric on ${\cal L}$ iff it satisfies the following
conditions:
\begin{itemize}

\vspace*{-2.25mm}\item[1)] $\lall\, a,b\in A,\, \xi ,\eta ,
\in{\cal L}\te\sb{Z(A)}A :\;
g(a\xi,\eta b)=ag(\xi,\eta)b$, where the left module structure on
${\cal L}\te\sb{Z(A)}A$ is given by $a(X\sb i\te\sb{Z(A)} c\sb i)
=X\sb i\te\sb{Z(A)} ac\sb i\,$.

\vspace*{-2.25mm}\item[2)] $\lall\, X,Y\in{\cal L} :\; g(X,Y)=g(Y,X)\,$
(symmetry condition).

\vspace*{-2.25mm}\item[3)] The induced map
$\widetilde{g}:{\cal L}\te\sb{Z(A)}A\ni\xi\longmapsto
g(\, .\, ,\xi )\in\Omega\sp 1(A)$ is an isomorphism of right $A$-modules.

\end{itemize}
\ede

\bde
A connection on ${\cal L}$ is said to be compatible with $g$ iff
\[
\lall\; X, Y ,Z\in{\cal L}\, :\;
Xg(Y,Z) = g(\nabla\!\sb XY ,Z) + g(Y,\nabla\!\sb XZ)\, .
\]
\ede

\bde\label{torsion}
The $Z(A)$-bilinear map
\[
T\sb\nabla :{\cal L}\times{\cal L}\ni(X,Y)\longmapsto
\nabla\sb XY-\nabla\sb YX-[X,Y]\ot\sb{Z(A)} 1\in{\cal L}\te\sb{Z(A)}A
\]
is called the torsion of $\nabla$.
\ede

\bde
Let $\,{\cal R}\sb\nabla(X,Y)$ be the
 $Z(A)$-homomorphism given by the formula:
\[
{\cal L}\ni Z\longmapsto{\cal R}\sb\nabla(X,Y)(Z)=
(\nabla\sp2 Y)(Z,X)\in{\cal L}\ot\sb{Z(A)}A\, ,
\]
and let
${\cal T}\sb E\in\mbox{\em Hom}\sb A\llp\mbox{\em Hom}\sb{Z(A)}
({\cal L},{\cal L}\te\sb{Z(A)}A),A\lrp$. We call the $Z(A)$-linear map
\[
Ric\sb\nabla:
{\cal L}\ni X\longmapsto{\cal T}\sb E\llp{\cal R}\sb\nabla(X,\, .\, )\lrp\in\hO\sp
1(A)
\]
the Ricci curvature of~$\nabla$.
\ede

\bre\label{projective}{\em
When ${\cal L}$ or $A$ is a finitely generated projective $Z(A)$-module,
\[
\mbox{Hom}\sb {Z(A)}({\cal L},{\cal L}\te\sb{Z(A)}A)
=\mbox{End}\sb {Z(A)}({\cal L})\ot\sb {Z(A)}A
\]
(see Proposition~2 in II.4\cite{bour}), and we can choose
${\cal T}\sb E$ to be a trace on $\mbox{End}\sb {Z(A)}({\cal L})$ 
tensored over
$Z(A)$ with $id\sb A$.
}\hfill{$\Diamond$}\ere

\bde
Let ${\cal M}({\cal L})$ and ${\cal C}({\cal L})$ denote the space of all
pseudo-Riemannian
metrics on~${\cal L}$, and the space of all connections on ${\cal L}$ respectively.
The functional $E:{\cal M}({\cal L})\times{\cal C}({\cal L})\lra k$ given by
the formula:
\[
E(g,\nabla)=-(\tau\sb{g}\ci{\cal T}\sb E)(\widetilde{g}\sp{-1}\ci Ric\sb\nabla)\,
,
\]
where $\tau\sb{g} : A\ra k$ is a metric dependent trace,
is called the Einstein action functional on~${\cal L}$.
\ede

\bre\em
With an appropriate choice of $\tau\sb g$ and ${\cal T}\sb E$
(cf.\ Proposition~\ref{main}), the functional $E$ coincides for
$A=C\sp\infty(M)$
and ${\cal L}=\mbox{Der}(C\sp\infty(M))$
with the standard Einstein--Hilbert action functional on $M$, for any
(paracompact) manifold $M$ admitting a 
(pseudo-)Riemannian metric.~\hfill{$\Diamond$}\ere

\section*{III. The case of $\boldmath Z(A)$-valued Metrics}

One can apply the same reasoning as in
the case of classical differential geometry to obtain:

\bpr[cf.~Section 9 in \cite{dm}]\label{lc}
Let $g({\cal L},{\cal L})\inc Z(A)$ or $\nabla\!\sb{\cal L}{\cal L}\inc{\cal L}$.
Then there exists at most one metric
compatible connection that is torsion free.
If it exists, it is given by the formula:\vspace*{-2.5mm}
\beq\label{lccf}
\nabla\sb XY=\frac{1}{2}\,\widetilde{g}\sp{-1}\lblp
Xg(Y,\, .\, )\! +\! Yg(X,\, .\, )\! -\! d\, g(X,Y)\!
 +\! g([X,Y],\, .\, )\! +\! g([\, .\, ,X],Y)\! +\! g([\, .\, ,Y],X)\lbrp\, .
\eeq
\epr\vspace*{2mm}

A connection given by (\ref{lccf}), will be called the Levi--Civita connection
of $g$, and denoted by $\nabla\!\sb g\,$.

For the rest of this section, we work under the assumption that 
$g({\cal L},{\cal L})\inc Z(A)$. 
The mathematical model considered here is practically identical
with a model constructed in \cite{akk,mkk}. (Compare (\ref{matrixact}) with
(3.24) in \cite{akk}, and Proposition~\ref{main} with (3.23), (3.21) in
\cite{akk} and (3.26) in~\cite{mkk}.) Note that there is an extra term in
(3.24)\cite{akk} that is absent in Proposition~\ref{main} due to our assumption
that
the classical and the algebraic derivations are orthogonal to each other.
Let us also mention that the dual point of view regarding linear connections
(i.e., where the space of 1-forms rather than ${\cal L}\te\sb{Z(A)}A$ is taken
as a starting point\footnote{
For duality issues of this kind see Section~6 in \cite{dm}.
}) was studied in \cite{mou} and \cite{mmm} also in the context of `matrix
geometry' (see Section~4.3 in \cite{mou} and Section~3 in~\cite{mmm}).

To begin with, let us consider the algebra of matrices $M\sb n({\Bbb R})$
and ${\cal L}=\mbox{Der}(M\sb n({\Bbb R}))={\frak s}{\frak l}(n,{\Bbb R})$.
Since ${\frak s}{\frak l}(n,{\Bbb R})$ is an $(n\sp2-1)$-dimensional vector space
over
${\Bbb R}$ (the centre of the algebra of matrices),
the endomorphisms of this space are
simply $(n\sp2-1)\times (n\sp2-1)$ matrices, and we can choose ${\cal T}\sb E$
to be the usual matrix trace tensored with $id\sb{M\sb n({\Bbb R})}$
(see Remark~\ref{projective}).

\bpr\label{matrix}
Let \mbox{
$A\! =\! M\sb n({\Bbb R}),\;{\cal L}\! =\! \mbox{\em Der}(M\sb n({\Bbb R})),\;
\tau\sb{g}\! =\! \frac{1}{n}\, |\det g|\sp{\frac{1}{2}}\, Tr$,
and~${\cal T}\sb E\! =\! Tr\te id\sb{M\sb n({\Bbb R})}$.} \\
Also, let $\{ E\sb i\}$ be any basis of         
${\frak s}{\frak l}(n,{\Bbb R})$.
Then
\beq\label{matrixact}
E(g,\nabla\!\sb g)=g\sp{jp}
(K\sb{jp}+\frac{1}{2}g\sp{il}g\sb{rk}c\sp r\sb{lp}c\sp k\sb{ij})\sqrt{|\det
g|}\, ,
\eeq
where $K$ is  the Killing metric on $SL(n,{\Bbb R})$, and
$\{K\sb{jp},g\sb{ij},g\sp{kl},c\sp r\sb{lp}\}$ are defined by the formulas
$K\sb{jp}=g(E\sb j,E\sb p),\; g\sb{rk}=g(E\sb r,E\sb k),\;
g\sp{pr}g\sb{rk}=\hd\sp p\sb k\, ,\; [E\sb l,E\sb p]=c\sp r\sb{lp}E\sb r$
respectively.
\epr

{\it Proof.} A direct computation using the symmetry of a metric and
the anti-symmetry of the Lie algebra structure constants.
\epf

\bre\em\label{flat}
Let $g$ be a left invariant metric on $SL(n,{\Bbb R})$. Then $g$ can be
identified
with a metric on  $\mbox{Der}(M\sb n({\Bbb R}))$. The connection on $SL(n,\Bbb
R)$
given by the left translations is compatible with all left invariant metrics.
The torsion part of this
connection is given by the formula (see (44) in \cite{tr}):
\[
{\frak T}\sp i\sb{jk}:=
\frac{1}{2}(Q\sp i\sb{jk}+g\sp{il}g\sb{jn}Q\sp n\sb{kl}
+g\sp{il}g\sb{kn}Q\sp n\sb{jl})
=\frac{1}{2}(c\sp i\sb{kj}+g\sp{il}g\sb{jn}c\sp n\sb{lk}
+g\sp{il}g\sb{kn}c\sp n\sb{lj})
\]
If the coefficients of the $M\sb n({\Bbb R})$--Levi--Civita connection of $g$
are defined by the equality\linebreak
\mbox{$\nabla\!\sb{E\sb j}E\sb k=\Gamma\sp i\sb{kj}E\sb i$}, then
\mbox{$\Gamma\sp i\sb{jk}={\frak T}\sp i\sb{jk}\in{\Bbb R}$}.
(Caution: One often defines the Christoffel symbols by the relation
\mbox{$\nabla\sb{E\sb j}E\sb k=\widetilde\Gamma\sp i\sb{jk}E\sb i$}.
In this notation, which is compatible with the notation used in \cite{tr},
the aforementioned relationship between the Christoffel symbols of the
noncommutative connection and the torsion part of the classical connection
can be equivalently written as
$\widetilde\Gamma\sp i\sb{jk}={\frak T}\sp i\sb{jk}+c\sp i\sb{jk}$.)
In general, if  Der$(A)$ equals the Lie algebra of some Lie group~$G$,
then the noncommutative torsion free \mbox{$g$-compatible} connection on
$\mbox{Der}(A)$ coincides in the above sense with the torsion part of the flat
connection on $G$ given by the left translations.
\hfill{$\Diamond$}
\ere

Our next step is to consider an algebra of matrix valued functions.
The module of derivations of such an algebra splits into two direct sum
components in the following way (cf.~Lemma~2.1 in \cite{dkm}):
\[
\mbox{Der}\llp C\sp{\infty}(M)\otimes M\sb n({\Bbb R})\lrp
=\mbox{Der}\llp C\sp{\infty}(M)\lrp\te{\Bbb R}\;\oplus\;
C\sp{\infty}(M)\te\mbox{Der}\llp M\sb n({\Bbb R})\lrp\, .
\]
For $M\! =\! T\sp m$, this module is a free
$C\sp{\infty}(T\sp m)$-module of dimension
$m+n\sp 2-1$. Consequently, its algebra of endomorphisms
is simply the algebra of matrices
$M\sb{m+n\sp 2-1}(C\sp{\infty}(T\sp m))$, and,
again, we can choose ${\cal T}\sb E$ to be the usual matrix trace (with values in
$C\sp{\infty}(T\sp m)\,$) tensored with $id\sb{M\sb n({\Bbb R})}$
(see Remark~\ref{projective}).

\bpr\label{main}
Let
$A=C\sp{\infty}(T\sp m)\otimes M\sb n({\Bbb R})\, ,\;{\cal L}=\mbox{\em Der}(A)
\, ,\;\tau\sb{g}=\frac{1}{n}\,\int\sb{T\sp m}|\det g|\sp{\frac{1}{2}}\;
Tr\,$, and ${\cal T}\sb E=Tr\ot id\sb{M\sb n({\Bbb R})}$.
Assume also that there exists a basis
$\{ E\sb i\}\sb{i\in\{ 1,\cdots ,m\}}$ of
$\,\mbox{\em Der}(C\sp{\infty}(T\sp m))$, and a basis
$\{ E\sb j\}\sb{j\in\{ m+1,\cdots m+n\sp 2-1\} }$ of
$\,\mbox{\em Der}(M\sb n({\Bbb R}))$
such that
\[
g\sb{ij}=\left\{
\begin{array}{ll}
0 & \mbox{for $i\leq m$ and $j>m$}\\
g\sb c(E\sb i,E\sb j) & \mbox{for $i,j\leq m$}\\
g\sb q(E\sb i,E\sb j) & \mbox{for $i,j>m$,}
\end{array}
\right.
\]
where $g\sb c$ is a classical (pseudo-)Riemannian metric on $T\sp m$, and $g\sb
q$
is a function that to each point of $T\sp m$ assigns a metric on
$\mbox{\em Der}(M\sb n({\Bbb R}))$.
Then
\beq\label{expr}
E(g,\nabla\!\sb g)=\int\sb{T\sp m}R\sb c\sqrt{|\det g\sb c|}\sqrt{|\det g\sb
q|}
\; +\int\sb{T\sp m}R\sb q\sqrt{|\det g\sb q|}\sqrt{|\det g\sb c|}
\; +\int\sb{T\sp m}R\sb{\mbox{\scriptsize\em mixed}}\sqrt{|\det g|}\, ,
\eeq
where $R\sb c$ is the classical scalar curvature of $g\sb c$, $R\sb q$
is a (point dependent) scalar curvature of $\nabla\!\sb{g\sb q}$ (i.e.\
\mbox{$R\sb q\! =\!{\cal T}\sb E(\widetilde{g}\sp{-1}\sb q\!\ci\!
Ric\sb{\nabla\!\sb{g\sb q}}\! )$}),
and $R\sb{\mbox{\scriptsize\em mixed}}$ is a function on $T\sp m$
that is a sum of mixed terms of the kind
\[
g\sp{AB}\frac{\partial}{\partial x\sp\mu}
\left( g\sp{\mu\nu}\frac{\partial g\sb{AB}}{\partial x\sp\nu}\right)\, ,
\;\; \mu ,\nu\leq m\, ,\;\; A,B > m\, .
\]
\epr

{\it Proof.} Very much like the proof of Proposition~\ref{matrix}.
\epf\ \\

As we see from Proposition~\ref{main}, the assumption that the metric $g$ is
\mbox{block-diagonal} allows us to split the Einstein action of the
Levi--Civita
connection into the following three terms:

\be

\item A \mbox{classical-like} term that differs from the usual
Einstein--Hilbert
 action on
$T\sp m$ only by the `quantum volume element' $\sqrt{|\det g\sb q|}$.

\item A quantum-like term that is equal to the integral over $T\sp m$
of the (point dependent) Einstein action of the Levi--Civita
connection on $M\sb n({\Bbb R})$.

\item A mixed term that involves the derivatives of  $g\sb q$.

\ee\vspace*{2mm}
{\parindent0pt
If  $g\sb q$ is constant over $T\sp m$, then $R\sb{\mbox{\scriptsize
mixed}}=0$,
and the expression (\ref{expr}) simplifies to
\[
E(g,\nabla\!\sb g)=
\sqrt{|\det g\sb q|}\; E(g\sb c)+E(g\sb q)\mbox{vol}(T\sp m,g\sb c)\, ,
\]
where $E(g\sb c)$ is the usual Einstein--Hilbert action,
$E(g\sb q)$ is the action computed
in Proposition~\ref{matrix}, and $\mbox{vol}(T\sp m,g\sb c)$ is the volume of
$T\sp m$ with respect to the metric $g\sb c$.
}

Finally, let us remark that the Yang--Mills (Maxwell) action calculated in
\cite{dkm}
for a similar algebra also splits into a classical-like, a quantum-like and
a mixed term.

\section*{IV. The Case of $\boldmath A$-valued Metrics}

Let us now lift the assumption that $g({\cal L},{\cal L})\inc Z(A)$, and
consider a toy model based on the following data:
\[
A=M\sb 4({\Bbb R})\, ,\;{\cal L}
={\frak s}{\frak o}(2)\oplus{\frak s}{\frak o}(2)\, ,\;
\tau\sb{g}=\mbox{${\displaystyle\frac{1}{4}}$}\,\sqrt{|\det g|}\;
Tr\sb{M\sb 4({\Bbb R})}\;\;\mbox{and}\;\;{\cal T}\sb E=Tr\ot id\sb{M\sb 4(\Bbb
R)}\, .
\]
We view ${\cal L}$ as a Lie subalgebra of $\mbox{Der}(M\sb 4({\Bbb R}))$
generated by $\hat F\sb 1:=[F\sb 1,\, .\, ]\, ,\;\hat F\sb 2:=[F\sb 2,\, .\,
]\, $,
where $F\sb 1:=\mbox{\scriptsize $\pmatrix{F&0\cr 0&0}$}\, ,\;
F\sb 2:=\mbox{\scriptsize $\pmatrix{0&0\cr 0&F}$}\, ,\;
F:=\mbox{\scriptsize $\pmatrix{0&1\cr -1&0}$}$.
A metric $g$ is treated as an element of
$GL\sb 2(M\sb 4({\Bbb R}))=GL\sb 8({\Bbb R})\,$. An assumption that $g$ is
symmetric reads
$g(\hat F\sb i,\hat F\sb j)=:g\sb{ij}=g\sb{ji}:=g(\hat F\sb j,\hat F\sb i)$,
or equivalently $g=g\sp T$,
where $\phantom{G}\sp T$ is the transpose in the algebra of $2\times 2$
matrices.
As to the inverse of $g$, denoted by
$g\sp{-1}=:
\mbox{\scriptsize $\pmatrix{g\sp{11}&g\sp{12}\cr g\sp{21}&g\sp{22}}$}\,$,
in general we do not have $g\sp{12}=g\sp{21}\,$. For instance, if
\[
g=\pmatrix
{2I\sb 2&0&I\sb 2&I\sb 2\cr 0&I\sb 2&I\sb 2&0\cr I\sb 2&I\sb 2&I\sb 2&0\cr
I\sb 2&0&0&I\sb 2}\, ,
\]
where $I\sb 2=\mbox{\scriptsize $\pmatrix{1&0\cr 0&1}$}\,$,
then $g\sp{12}\not= g\sp{21}\,$. Very much as we did before, we define the
Christoffel symbols by
$\nabla\!\sb{\hat F\sb i}\hat F\sb j=\hat F\sb k\te\hG\sp k\sb{ji}\,$,
the curvature coefficients by
$(\nabla\sp 2\hat F\sb k)(\hat F\sb i,\hat F\sb j)
=\hat F\sb m\te R\sp m\sb{kij}\,$,
and the Ricci curvature coefficients by
$(Ric\sb\nabla\hat F\sb i)(\hat F\sb j)=R\sb{ji}\,$.
It is straightforward to verify that
$R\sp m\sb{kij}=\hG\sp m\sb{ni}\hG\sp n\sb{kj}-\hG\sp m\sb{nj}\hG\sp n\sb{ki}+
[F\sb i,\hG\sp m\sb{kj}]-[F\sb j,\hG\sp m\sb{ki}]$
(see Proposition~\ref{spera}), $R\sb{kj}=R\sp i\sb{kij}\,$, and
\beq\label{act4}
E(g,\nabla)=-\frac{1}{4}\sqrt{|\det g|}\; Tr\sb{M\sb 4(\Bbb
R)}(g\sp{jk}R\sb{kj})
=-\frac{1}{4}\sqrt{|\det g|}\; tr(g\sp{-1}r)\, ,
\eeq
where $tr$ denotes the usual trace on $M\sb 8({\Bbb R})$, and
$r:=\mbox{\scriptsize $\pmatrix{R\sb{11}&R\sb{12}\cr R\sb{21}&R\sb{22}}$}\,$.
In what follows, rather than look for a Levi--Civita connection,
we will use an analogue of the Palatini method of variation
 (see 21.2 in \cite{mtw}, cf.~5.4--5.5 in \cite{gs})
and find critical points of the Einstein action functional.
First, let us determine
the equation equivalent to the condition that the variation of $E$ with respect
to $g$ vanish. At any critical point, for an arbitrary matrix
$h\in M\sb 8({\Bbb R})$ with the coefficients
satisfying $h\sb{12}=h\sb{21}\in M\sb 4({\Bbb R})$ we must have:
\beq\label{ddsg}
\frac{d}{ds}E(g+sh,\nabla)|\sb{s=0}=0\, .
\eeq
After substituting (\ref{act4}) into (\ref{ddsg}), and carrying out
the differentiation, we obtain:
\beq\label{f1}
\frac{1}{2}tr(g\sp{-1}r)tr(hg\sp{-1})=tr(hg\sp{-1}rg\sp{-1})\, .
\eeq
Now, since (\ref{f1}) must be true for any matrix $h$ such that
$h\sb{12}=h\sb{21}$, we can conclude that
\[
\lall\, i,j\in\{ 1,2\}:\;\; \frac{1}{2}tr(g\sp{-1}r)(g\sp{ij}+g\sp{ji})=
g\sp{ik}R\sb{kl}g\sp{lj}+g\sp{jk}R\sb{kl}g\sp{li}\, .
\]
Multiplying both sides by $g\sb{im}$, and then taking the trace, we find that
\[
8tr(g\sp{-1}r)=2tr(g\sp{-1}r)\, .
\]
Hence, very much as in the classical general relativity, we have
$tr(g\sp{-1}r)=0$. Consequently,
\beq\label{field}
g\sp{-1}rg\sp{-1}+(g\sp{-1}rg\sp{-1})\sp T =0\, .
\eeq
The formula (\ref{field}) is an analogue of the Einstein vacuum field equation.
(Observe that the implication $(\ref{field})\Rightarrow (\ref{ddsg})$ is also
true.)

Our next step is to find an explicit form of the equations equivalent to the
vanishing of the variation of $E$ with respect to $\nabla$. Let $\nabla + sA$
denote a connection on ${\cal L}$ whose Christoffel symbols are
$\hG\sp{k}\sb{ji}+sA\sp{k}\sb{ji}\,$, where
$s\in{\Bbb R}\, ,\; A\sp{k}\sb{ji}\in M\sb 4({\Bbb R})\, ,\; i,j,k\in \{ 1,2\}$.
Then the condition that
\[
\frac{d}{ds}E(g,\nabla +sA)|\sb{s=0}=0\;\; \mbox{for any}\;\; A
\]
is equivalent to the following 8 equations:

\beq\label{c1}
g\sp{11}\hG\sp{2}\sb{12}+\hG\sp{1}\sb{22}g\sp{22}+
[\hG\sp{1}\sb{12}+F\sb 2,g\sp{21}]
=0\eeq\beq\label{c2}
g\sp{11}\hG\sp{2}\sb{11}+\hG\sp{1}\sb{21}g\sp{22}+
[\hG\sp{1}\sb{11}+F\sb 1,g\sp{21}]
=0\eeq\beq\label{c3}
g\sp{22}\hG\sp{1}\sb{22}+\hG\sp{2}\sb{12}g\sp{11}+
[\hG\sp{2}\sb{22}+F\sb 2,g\sp{12}]
=0\eeq\beq\label{c4}
g\sp{22}\hG\sp{1}\sb{21}+\hG\sp{2}\sb{11}g\sp{11}+
[\hG\sp{2}\sb{21}+F\sb 1,g\sp{12}]
=0\eeq\beq\label{c5}
g\sp{22}\hG\sp{1}\sb{12}-\hG\sp{2}\sb{22}g\sp{22}-
g\sp{12}\hG\sp{2}\sb{12}-\hG\sp{2}\sb{12}g\sp{21}-
[F\sb 2,g\sp{22}]=0\eeq\beq\label{c6}
g\sp{11}\hG\sp{2}\sb{22}-\hG\sp{1}\sb{12}g\sp{11}-
g\sp{21}\hG\sp{1}\sb{22}-\hG\sp{1}\sb{22}g\sp{12}-
[F\sb 2,g\sp{11}]=0\eeq\beq\label{c7}
g\sp{22}\hG\sp{1}\sb{11}-\hG\sp{2}\sb{21}g\sp{22}-
g\sp{12}\hG\sp{2}\sb{11}-\hG\sp{2}\sb{11}g\sp{21}-
[F\sb 1,g\sp{22}]=0\eeq\beq\label{c8}
g\sp{11}\hG\sp{2}\sb{21}-\hG\sp{1}\sb{11}g\sp{11}-
g\sp{21}\hG\sp{1}\sb{21}-\hG\sp{1}\sb{21}g\sp{12}-
[F\sb 1,g\sp{11}]=0
\eeq

It is straightforward to verify that a connection $\nabla\sp g$
(not to be confused with the Levi--Civita connection $\nabla\!\sb g$)
whose Christoffel symbols are given by

\[
\hG\sp{1}\sb{22}=\hG\sp{1}\sb{21}=-g\sp{11}\;\;\;\;\;\;\;\;\;\;\;\;\;\;\;\;
\hG\sp{2}\sb{22}=-F\sb 2-g\sp{12}\;\;\;\;\;\;\;\;\;\;\;\;\;\;\;\;
\hG\sp{1}\sb{12}=-F\sb 2+g\sp{21}
\]\[
\hG\sp{2}\sb{12}=\hG\sp{2}\sb{11}=\phantom{-}g\sp{22}\;\;\;\;\;\;\;\;\;\;\;\;\;\;\;\;
\hG\sp{2}\sb{21}=-F\sb 1-g\sp{12}\;\;\;\;\;\;\;\;\;\;\;\;\;\;\;\;
\hG\sp{1}\sb{11}=-F\sb 1+g\sp{21}
\]\vspace*{.01mm}

satisfies (\ref{c1}--\ref{c8}), and has vanishing Ricci curvature
($Ric\sb{\nabla\sp g}=0$).
We have thus arrived at the following:
\bpr\label{cpe}
Let $A,\, {\cal L},\, \tau\sb g,\,$ and ${\cal T}\sb E$ be as above. Then
\[
\lall\, g\in{\cal M}({\cal L})\;\lex\,\nabla\sp g\in{\cal C}({\cal L})
\mbox{such that}\;\;
(g,\nabla\sp g)\;\;\mbox{is a critical point of}\;\; E.
\]
Furthermore, the value of $E$ at any critical point is zero.
\epr

\bre\em
The value of the functional $E$ calculated at non-critical points is not
necessarily zero. For example, take a metric $g\sb 0=g\sb 0\sp{-1}$
with the components
\[
(g\sb 0)\sb{11}=(g\sb 0)\sb{22}=0\, ,\;\;
(g\sb 0)\sb{12}=(g\sb 0)\sb{21}=\pmatrix{I\sb 2&0\cr 0&K}\, ,\;\;
K:=\pmatrix{0&1\cr 1&0}\, ,
\]
 and take a connection
$\nabla\sb 0$ whose only non-vanishing Christoffel symbol is
$(\hG\sb 0)\sp{1}\sb{11}=\mbox{\scriptsize $\pmatrix{0&0\cr 0&J}$}\,$, where
$J=\mbox{\scriptsize $\pmatrix{1&0\cr 0&-1}$}$. Then $E(g\sb 0,\nabla\sb
0)=-1$.
\hfill{$\Diamond$}\ere

\bre\em
Note that  for $\tau\sb g$ equal to $Tr\sb{M\sb 4({\Bbb R})}$ rather than to
$\frac{1}{4}\, |\det g|\sp{\frac{1}{2}}\, Tr\sb{M\sb 4({\Bbb R})}$
the field equation (\ref{field}) and the equations (\ref{c1}--\ref{c8}) are
still
satisfied. Also, $tr(g\sp{-1}r)$ still equals zero at any critical point.
Consequently, Proposition~\ref{cpe} remains valid as well, if we replace
$\frac{1}{4}\, |\det g|\sp{\frac{1}{2}}\, Tr\sb{M\sb 4({\Bbb R})}$ by the usual
trace on $M\sb 4({\Bbb R})$.
\hfill{$\Diamond$}\ere

\bre\em
The fact that the functional $\widetilde E:{\cal M}({\cal L})\ni g\longmapsto
E(g,\nabla\sp g)\in {\Bbb R}$ equals identically zero is a reflection of the
same effect that we observe for the usual 2-torus. We might try to push this
analogy even further and say that we think of a circle $S\sp 1$  as a Lie group
generated by ${\frak s}{\frak o}(2)$, and replace $C\sp{\infty}(S\sp 1)$ by
$M\sb 2({\Bbb R})$ for which ${\frak s}{\frak o}(2)$ is the space of all derivations
satisfying $X(a\sp T)=(Xa)\sp T$. Then it is natural to replace
$C\sp{\infty}(T\sp 2)$ by $M\sb 2({\Bbb R})\te M\sb 2({\Bbb R})=M\sb 4({\Bbb R})$.
\hfill{$\Diamond$}\ere

Observe that although $\nabla\sp g$ functions as if it were a
Levi--Civita connection, it is in general neither metric nor torsion free
(e.g., take $g$ to be the identity matrix of $GL\sb 8({\Bbb R})\,$).
It is perhaps worth emphasizing, however, that the metric compatibility
condition,
which can be equivalently written as:
\[
g\sp{pj}\hG\sp{n}\sb{ji}+\hG\sp{p}\sb{ji}g\sp{jn}+[F\sb i,g\sp{pn}]=0\, ,
\;\; i,p,n\in\{ 1,2\}\, ,
\]
is not very different from the formulas (\ref{c1}-\ref{c8}).
It would be interesting to
find a functional on ${\cal M}({\cal L})\times{\cal C}({\cal L})$ that not only would
coincide with the usual Einstein--Hilbert functional 
(or whose equations of motion would agree with the standard ones)
in the case of classical
geometry, but also would yield, through its variation with respect to
connection,
the metric compatibility condition.
}\\

\small
{\bf Acknowledgments:} 
This work was supported in part by a visiting fellowship at
the International Centre for Theoretical Physics in Trieste,
and the KBN grant 2 P301 020 07. It is my great pleasure to thank 
Andrzej Borowiec, Marc Rieffel and Tsou Sheung Tsun
for their manifold assistance and encouragement, Michel Dubois-Violette for an
electronic mail discussion, and Giovanni Landi and Andrzej Sitarz for their
help in proofreading the manuscript.

\end{document}